\documentclass[12pt]{article}

\usepackage{german}
\usepackage{amsmath}
\usepackage{amssymb}
\usepackage{latexsym}

\newcommand{\ang}{\raisebox{0.2ex}{\scriptsize$\triangleright$}}

\newcommand{\mn}{\medskip\noindent}
\newcommand{\bn}{\bigskip\noindent}
\newcommand{\sn}{\smallskip\noindent}

\newcommand{\rti}{\,{\scriptstyle\rtimes}\,} %

\newcommand{\D}{{\mathcal{D}}} 
\newcommand{\E}{{\mathcal{E}}}

\newcommand{\Hh}{{\mathcal{H}}}
\newcommand{\cO}{{\mathcal{O}}}
\newcommand{\K}{{\mathcal{K}}}
\newcommand{\U}{{\mathcal{U}}}
\newcommand{\T}{{\mathcal{T}}}
\newcommand{\X}{{\mathcal{X}}}
\newcommand{\A}{{\mathcal{A}}}

\newcommand{\B}{{\mathcal{B}}}

\newcommand{\N}{\mathbb{N}}
\newcommand{\R}{\mathbb{R}}
\newcommand{\C}{\mathbb{C}}

\newcommand{\Lin}{{\mathrm{Lin}}}
\newcommand{\id}{{\mathrm{id}}}

\newcommand{\br}{{\mathbf{r}}}

\begin{document}
\vspace{10.6mm}
\title{}\date{}
\begin{center}
{\bf {\Large{Commutator Representations of Covariant Differential Calculi on Quantum Groups}}}

\bn
{\large Konrad Schm"udgen}

\sn
\small{Fakult"at f"ur Mathematik und Informatik\\ Universit"at Leipzig, Augustusplatz 10, 04109 Leipzig, Germany\\ E-mail: schmuedg@mathematik.uni-leipzig.de }
\end{center}

\sn
{\bf Abstract:} Let $(\Gamma,d)$ be a first order differential 
$\ast$-calculus on a $\ast$-algebra $\A$. We say that a pair $(\pi,F)$ of a $\ast$-representation $\pi$ of $\A$ on a dense domain $\D$ of a Hilbert space 
and a symmetric operator $F$ on $\D$ gives a commutator representation of 
$\Gamma$ if there exists a linear mapping $\tau:\Gamma \to L(\D)$ 
such that $\tau (adb)=\pi(a)i[F,\pi(b)]$, $a,b\in\A$. Among others, 
it is shown that each left-covariant $\ast$-calculus $\Gamma$ of a compact 
quantum group Hopf $\ast$-algebra $\A$ has a faithful commutator 
representation. 
For a class of bicovariant $\ast$-calculi on $\A$ there is a commutator 
representation such that $F$ is the image of a central element of the quantum 
tangent space. If $\A$ is the Hopf $\ast$-algebra of the compact form of 
one of the quantum groups $SL_q(n{+}1)$, $O_q(n)$, $Sp_q(2n)$ with real 
trancendental $q$, then this commutator representation is faithful. 

\sn

\mn
{\bf 0. Introduction}\\ In the last decade covariant differential calculi on 
quantum groups have been investigated and a general theory of such calculi 
has been developed ([W], see [KS], Chapter 14, for a 
thorough treatment). Though this theory fits nicely into the general 
framework of quantum groups, its relations to Alain Connes' 
noncommutative geometry [C] are not yet understood. In [S1] we have 
shown that the distinguished covariant differential calculi on the 
quantum group $SU_q(2)$ cannot be described by a spectral triple. However they can be given by means of {\it unbounded} commutators 
in Hilbert space representations. 

In this paper we address this problem in a more general setting. To be 
precise, suppose that $\Gamma$ is a first order differential 
$\ast$-calculus on a $\ast$-algebra $\A$ with differentiation
 $d:\A\rightarrow \Gamma$. Let $\pi$ be a $\ast$-representation of 
$\A$ by possibly unbounded operators on a dense invariant domain 
$\D$ of a Hilbert space and let $F$ be a symmetric operators 
$\D$ which leaves $\D$ invariant. If there exists a well-defined 
linear mapping $\tau$ of $\Gamma$ into the linear operators on $\D$ such that 
\begin{equation}\label{tau}
\tau (adb)=\pi(a)i[F,\pi(b)],a,b\in\A,
\end{equation}
then the pair $(\pi,F)$ is called a {\it commutator representation} 
of the $\ast$-calculus $\Gamma$. We also study the problem at purely 
algebraic level: If $\pi$ is a homomorphism of $\A$ into another 
algebra $\B$ and $F$ is an element of $\B$ such that there exist a 
linear map $\tau$: $\A\rightarrow \B$ satisfying (\ref{tau}), then 
we say that the pair $(\pi,F)$ is an {\it algebraic commutator 
representation} of $\Gamma$.

We briefly state some of our main results. Each 
left-covariant first order calculus on a Hopf algebra $\A$ admits a 
faithful algebraic commutator representation (Proposition 1). A large 
class of bicovariant calculi on coquasitriangular Hopf algebras $\A$ has 
an algebraic commutator representation by a single element of the quantum 
tangent space which in central in the Hopf dual $\A^\circ$ (Proposition 4). 
For the coordinate Hopf algebras $\A=\cO(G_q), G_q=SL_q(n{+}1), O_q(n), 
Sp_q(2n)$, with transcendental $q$ this representation is faithful (Theorem 5). 

A crucial role in these considerations plays the cross product 
algebra $\A\rti\A^\circ$. Suppose that $\A$ is a 
Hopf $\ast$-algebra. Combined with $\ast$-representations of the 
$\ast$-algebra $\A\rti\A^\circ$ (more precisely, of an appropriate 
$\ast$-subalgebra $\U$) the above algebraic commutator representations lead 
to Hilbert space commutator representations (Theorems 2 and 7). For 
instance, if $\A$ is a $CQG$ algebra (compact quantum group algebra), 
we can take the Heisenberg representation of the cross product 
$\ast$-algebra $\A\rti\A^\circ$. In this manner we obtain a 
faithful commutator representation of each left-covariant 
$\ast$-calculus (Corollary 3) and of bicovariant $\ast$-calculi 
on the quantum groups $SL_q(n{+}1)$, $O_q(n)$ and $Sp_q(2n)$ for real 
transcendental $q$ (Corollary 8). 

In the final Section 4 we list three simple examples of commutator 
representations of covariant differential calculi on quantum spaces.

Since the FODC $\Gamma$ is left-covariant with respect to $\A$, it is right-covariant with respect to $\U$. For a commutator 
representation of $\Gamma$ one would like to have  
the $\U$-symmetry on the algebra $L(\D)$ as well. For the 
commutator representations constructed by cross product algebras 
this can be done, while for the examples in Section 4  
representations in larger Hilbert spaces are needed. These problems will 
be studied in a forthcoming paper with E. Wagner. 

The results of the paper have been presented at workshops in 
Bayrischzell (March 2001) and Berkeley (April 2001). 

\sn
{\bf Acknowledgement:} I am indebted to I. Heckenberger, L. Vaksman and E. Wagner for 
instructive discussions on the subject of the paper.
 
\bn
{\bf 1. Preliminaries}\\
In this paper all algebras are unital and over the complex field. Algebra 
homomorphisms are always unit preserving. Suppose $\A$ is an algebra. We 
denote by $M_n(\A)$ the algebra of $n{\times}n$-matrices with entries 
from $\A$ and by $\A_n$ the direct sum algebra 
$\A\oplus{\dots}\oplus \A$ ($n$ times). The algebra of linear mappings 
of a vector space $V$ is denoted by $L(V)$. 
A {\it first order differential calculus} (abbreviated, a FODC) over 
$\A$ is an $\A$-bimodule $\Gamma$ equipped with a linear mapping 
$d:\A\rightarrow\Gamma$ such that $d(ab)=a{\cdot}db+da{\cdot}b$ 
for $a,b\in\A$ and $\Gamma=\Lin\{a{\cdot}db; a,b\in\A\}$.

Let $\rho:\A\rightarrow B$ be an algebra homomorphism of $\A$ into another 
algebra $\B$ and let $C\in\B$. We shall say that the pair $(\rho,C)$ is 
an {\it algebraic commutator representation} of a FODC $\Gamma$ over $\A$ 
if there exists a well-defined (!) linear mapping 
$\tau:\Gamma\rightarrow\B$ such that
\begin{equation}\label{acom}
\tau(adb)=i\rho(a)(C\rho(b)-\rho(b)C), a,b\in\A.
\end{equation}
If $\tau$ is injective, then the pair 
$(\rho,C)$ is called {\it faithful}. (The complex unit $i$ in (\ref{acom}) 
is only included for a convenient treatment of $\ast$-calculi.)

Now suppose that $\A$ is a $\ast$-algebra. A FODC over $\A$ is called a 
$\ast$-{\it calculus} if there exist a well-defined (!) involution 
$\gamma\rightarrow\gamma^\ast$ of the vector space $\Gamma$ such 
that $(adbc)^\ast=c^\ast\cdot d(b^\ast)\cdot a^\ast$ for all $a,b,c\in\A$. 
A {\it $\ast$-representation} of $\A$ is an algebra homomorphism $\pi$ of 
$\A$ into the algebra $L(\D)$ of linear operators of a dense linear 
subspace $\D$ of a Hilbert space such that
$$
\langle\pi(a)\varphi,\psi\rangle=\langle\varphi,\pi(a^\ast)\psi\rangle, 
\varphi,\psi\in\D,a\in\A.
$$
We say that a pair $(\pi,F)$ of a $\ast$-representation $\pi$ of $\A$ on 
the domain $\D$ and a symmetric linear operator $F\in L
(\D)$ is a {\it commutator representation} of a $\ast$-calculus $\Gamma$ 
over $\A$ if $(\pi,F)$ is an algebraic commutator representation of 
$\Gamma$ with $\B=L(\D)$.

Suppose  that $\A$ and $\U$ are Hopf algebras and $\langle\cdot,\cdot\rangle:\U\times\A\rightarrow\C$ is a dual 
pairing of Hopf algebras. Then $\A$ is a left $\U$-module algebra 
with left action $\ang$ of $\U$ on $\A$ given by
$$
f\ang a=a_{(1)}\langle f,a_{(2)}\rangle, f\in\U, a\in\A,
$$
where $\Delta(a)=a_{(1)}\otimes a_{(2)}$ is the Sweedler notation of 
the comultiplication of $\A$. The {\it left cross product algebra} 
$\A\rti\U$ is the vector space $\A\otimes \U$ equipped with the product 
defined by
\begin{equation*}
(b\otimes f)(a\otimes g)=b(f_{(1)}\ang a)\otimes f_{(2)}g=\langle f_{(1)}, a_{(2)}\rangle b a_{(1)}\otimes f_{(2)} g
\end{equation*}
for  $a,b\in\A$ and $f,g\in\U$, or equivalently $\A\rti\U$ is the algebra 
generated by the two subalgebras $\A$ and $\U$ with cross relation
\begin{align}\label{left}
fa=(f_{(1)}\ang a)f_{(2)}=\langle f_{(1)}, a_{(2)}\rangle a_{(1)} f_{(2)},~ 
a\in\A, f\in \U.
\end{align}

\bn
{\bf 2.Commutator representation of left-covariant FODC on Hopf algebras}\\
In this section $\Gamma$ is a finite dimensional left-covariant FODC of a Hopf algebra $\A$. We freely use some facts on covariant FODC; see [KS], 14.1.

Let $\{\omega_1,{\dots},\omega_n\}$ be a basis of the vector space of 
left-invariant elements of $\Gamma$. Then there exists functionals 
$X_k, f^k_l, k,l=1,{\dots},n$, of the dual Hopf algebra $\A^\circ$ of 
$\A$ such that 

\begin{align}\label{da}
&da=\sum_k (X_k\ang a) \omega_k,\\
\label{omegaa}
&\omega_k a=\sum_j(f^k_j\ang a)\omega_j,\\
\label{comult}
&\Delta(X_k)=\varepsilon\otimes X_k+\sum_j X_j\otimes f^j_k
\end{align}
for $a\in \A$. The linear span $\T_\Gamma$ of functionals 
$X_1,{\dots},X_n$ is the quantum tangent space of the left-covariant 
FODC $\Gamma$. From (\ref{comult}) and (\ref{left}) it follows that the 
functionals $X_k, f^k_j$ and elements $a\in\A$ satisfy the commutation 
relations
\begin{align}\label{xfrel}
&X_ka=aX_k+\sum_j(X_j\ang a) f^j_k\\
\label{ffrel}
&f^k_j a=\sum_l(f^k_l\ang a) f^l_j
\end{align}
in the cross product algebra $\A\rti\A^\circ$. There is a left action $.$ of the algebra $\A\rti\A^\circ$ on $\A$ given by $f .  b=f\ang b$ and $a . b=a{\cdot} b$ for $f\in\A^\circ$ and $a,b\in\A$. By matrix multiplication this yields a left action of the matrix algebra $M_{n+1}(\A\rti\A^\circ)$ on $\A_{n+1}=\A\oplus{\dots}\oplus\A$. In this manner, $M_{n+1}(\A\rti\A^\circ)$ becomes a subalgebra of $L (\A_{n+1})$. Define elements $C$ and $\Omega_k,k=1,{\dots},n$, of $M_{n+1}(\A^\circ)\subseteq M_{n+1} (\A\rti\A^\circ)\subseteq L(\A_{n+1})$ by
\begin{align*}
C=\left( \begin{matrix} 0      &X_1    &{\cdots}&X_n\\
                        X_1    &0      &{\cdots}&0\\
                        \vdots &\vdots &        &\vdots\\
                        X_n    &0      &{\cdots}&0
         \end{matrix} \right)~,~~
\Omega_k=i\left( \begin{matrix}  
                        0     &f^k_1 &{\cdots} &f^k_n\\
                        f^k_1 &0     &{\cdots} &0\\
                        \vdots&\vdots&         &\vdots\\
                        f^k_n &0     &{\cdots}  &0
          \end{matrix}\right)
\end{align*}
Using equations (\ref{da}), (\ref{omegaa}), (\ref{xfrel}) and (\ref{ffrel}) we derive that 
\begin{align}\label{r1}
&i(C a-aC) . (b_0,{\dots},b_n)=\sum_k(X_k\ang a)\Omega_k . (b_0,{\dots}, b_n),\\
\label{r2}
&\Omega_k a . (b_0,{\dots},b_n)=\sum_l(f^k_l\ang a)\Omega_l . (b_0,{\dots}, b_n)
\end{align}
for $a,b_0,{\dots},b_n\in\A$ and $k=1,{\dots},n$. Define $\varphi (a)(b_0,{\dots},b_n)=(ab_0,{\dots},ab_n)$. From the relations (\ref{r1}) and (\ref{r2}) we obtain

\mn
{\bf Proposition 1.}  {\it The pair $(\rho,C)$ is a faithful algebraic commutator representation of the left-covariant FODC $\Gamma$ over $\A$ with $\B=L(\A_{n+1})$.}

Now we turn to Hilbert space commutator representations. We suppose that 
$\A$ is a Hopf $\ast$-algebra and $\Gamma$ is a $\ast$-calculus.

Then the quantum tangent space $\T_\Gamma$ of $\Gamma$ is $\ast$-invariant ([KS], Proposition 14.6). Hence we can assume without loss of generality that all functionals $X_i$ are hermitean, that is, $X^\ast_j=X_j$ for $j=1,{\dots},n$. Further, since $\A$  is a left module $\ast$-algebra for the Hopf algebra $\A^\circ$, the cross product algebra $\A\rti\A^\circ$ is a $\ast$-algebra which contains $\A$ and $\A^\circ$ as $\ast$-subalgebras ([SW], Lemma 2.1). Let $\U_0$ be a $\ast$-subalgebra of $\A^\circ$ which contains all functionals $X_j,f^k_j,k,j=1,{\dots},n$. For notational simplicity we denote the $\ast$-subalgebra of $\A\rti \A^\circ$ generated by $\A$ and $\U_0$ by $\A\rti\U_0$.

Suppose that $\pi_1$ is a $\ast$-representation of the $\ast$-algebra 
$\A\rti\U_0$ on the domain $\D$. It yields a $\ast$-representation 
$\pi_{n+1}$ of the $\ast$-algebra $M_{n+1}(\A\rti\U_0)$ on the domain 
$\D_{n+1}:=\D\oplus{\dots}\oplus\D(n{+}1~ {\rm times})$ defined by 
$\pi_{n+1}((x_{kj})):=((\pi(x_{kj})))$ for $(x_{kj})\in M_{n+1} (\A\rti\
U_0)$. 
Let $\pi$ denote the $\ast$-representation of $\A$ on $\D_{n+1}$ 
given by $\pi(a)(\varphi_0,{\dots},\varphi_n) := (\pi_0(a)\varphi_0,
{\dots}, \pi_0(a)\varphi_n),\varphi_0,{\dots}, \varphi_n\in\D$. 
Since $X^\ast_j=X_j$ for $j=1,{\dots},n$, the operator $F:=\pi_{n+1}(C)$ on the domain $\D_{n+1}$ is symmetric. Since $\pi_{n+1}$ is an 
algebra homomorphism of $M_{n+1}(\A\rti\U^0)$ into $L(\D_{n+1})$, relations (\ref{r1}) and (\ref{r2}) remain valid for the corresponding images in Hilbert space. That is, we have proved

\sn
{\bf Theorem 2.} {\it Suppose $\Gamma$ is a finite dimensional left-covariant first order $\ast$-calculus of the Hopf $\ast$-algebra $\A$. If $\pi_1$ is a $\ast$-representation of the $\ast$-algebra $\A\rti\U_0$ defined above, then the pair $(\pi,F)$ is a commutator representation of $\Gamma$. It is faithful when the $\ast$-representation $\pi_1$  is faithful.}

By the preceding theorem, $\ast$-representations of cross product 
algebras allow us to construct commutator representations of left-covariant 
FODC. This observation was one of the key motivations for the 
study of Hilbert space representations of cross product algebras in [SW].

Now suppose that $\A$ is a $CQG$ algebra, that is, $\A$ is a Hopf $\ast$-algebra which is the linear span of matrix elements of finite dimensional unitary corepresentations of $\A$ (see e.g. [KS], 11.3). Let $h$ denote the Haar state of $\A$ and let $\pi_h$ denote the Heisenberg representation of the cross 
product $\ast$-algebra $\A\rti\A^\circ$ ([SW], 5.2). That is, $\pi_h$ is the unique $\ast$-representation of $\A\rti\A^\circ$ such that its restriction to $\A$ is the GNS representation of the state $h$ on $\pi_h(\A)\varphi_n$ and $\pi_h(f)\varphi_h=f(1)\varphi_h$ for $f\in\A^\circ$. Applying Theorem 2 with $\pi_1=\pi_h\lceil\A\rti\U_0$ we obtain

\mn
{\bf Corollary 3.} {\it Each finite dimensional left-covariant first order $\ast$-calculus of a CQG algebra has a faithful commutator representation.}

\bn
{\bf 3. Commutator representation of bicovariant FODC on coquasi\-triangular Hopf algebras}

In this section we suppose that $\A$ is a coquasitriangular Hopf algebra 
with universal $\br$-form $\br$ (see e.g. [KS], Chapter 10, for this notion). Let $\bar{\br}$ denote the convolution inverse of $\br$.
 First let us briefly repeat the general 
construction of bicovariant FODC over $\A$ developped in [KS], 14.5--14.6. 
In a special case this method was invented in [Ju].

Fix a matrix corepresentation $v{=}(v^k_l)_{k,l=1,{\dots},n}$ and a central character $\zeta$ of $\A$. That is, $v{=}(v^k_l)$ is a matrix from $M_n(\A)$ satisfying $\Delta(v^k_l){=}\sum_jv^k_j\otimes v^j_l$ and $\varepsilon(v^k_l){=}\delta_{kl}, k,l{=}1,{\dots},n$ and $\zeta$ is an algebra homomorphism of $\A$ into $\C$. 

Define functionals  $l^{\pm k}_j, l^k_j\in\A^\circ$ by
$$
l^{+k}_j(\cdot)=\br(\cdot\otimes v^k_j),l^{-k}_j(\cdot)=\bar{\br}(v^k_j\otimes \cdot), l^k_j=\sum_tS(l^{-k}_t)l^{+t}_j.
$$
Let $L^+,L^{-,c}\in M_n(\A^\circ)$ and $v^c\in M_n(\A)$ be the matrices 
with $(k,j)$ entries $l^{+k}_j,S(l^{-j}_k)$ and $S(v^j_k)$, respectively. Set
\begin{equation}\label{fs}
f^{kj}_{st}:=\zeta S(l^{-s}_k)l^{+j}_t, k,j,s,t=1,{\dots},n.
\end{equation}
Then, $\tilde{\Gamma}_{v,\zeta}:=(1\otimes v\otimes v^c,\zeta\otimes 
L^{-,c}\otimes L^+)$ is a Hopf bimodule of $\A$. On the canonical basis $\{\theta_{kj}, k,j=1,{\dots},n\}$ of left-invariant elements the $\A$-bimodule structure of $\tilde{\Gamma}_{v,\zeta}$ is determined by the relations
\begin{equation}\label{comm1}
\theta_{kj} a=\sum_{s,t} (f^{kj}_{st}\ang a)\theta_{st}, a\in\A.
\end{equation}
The element $\theta:=\sum_k\theta_{kk}$ is biinvariant. Defining 
\begin{equation}\label{diff}
da:=\theta a-a\theta,a\in\A,
\end{equation}
$\Gamma_{v,\zeta}:=\A\cdot d\A\cdot\A$ is a bicovariant FODC over $\A$ with differentiation map $d$. From (\ref{fs}) and (\ref{comm1}) we conclude that the quantum tangent space of $\Gamma_{v,\zeta}$ is the linear span of functionals
\begin{equation}\label{xkj}
X_{kj}=\zeta l^k_j-\delta_{kj}\varepsilon, k,j=1,{\dots},n.
\end{equation}
Put $A^j_k:=\sum_s\br(S^2(v^j_s)\otimes v^s_k)$. By Proposition 4.16(ii) in [KS], the functional
\begin{equation}\label{cv}
C_{v,\zeta}:=\sum_{k,j}X_{kj} A^j_k
\end{equation}
of the quantum tangent space of $\Gamma_{v,\zeta}$ belongs to the center of the dual Hopf algebra $\A^\circ$.

Now we construct an algebraic commutator representation of the FODC $\Gamma_{v,\zeta}$ in the cross product algebra $\A\rti\A^\circ$. 
Set $\Omega_{kj}=\sum_{s,t}\zeta S(l^{-s}_k)l^{+j}_t A^t_s$.
Since 
$$
\Delta(\Omega_{kj})=\sum_{i,l,s,t}\zeta S(l^{-i}_k)l^{+j}_l\otimes \zeta S(l^{-s}_i)l^{+l}_t A^t_s=\sum_{i,l} f^{kj}_{il} \otimes \Omega_{il},
$$ 
it follows from (\ref{left}) that the elements $\Omega_{kj}$ and $a\in\A$ satisfy the commutation relations
\begin{equation}\label{comm2}
\Omega_{kj} a=\sum_{i,l}(f^{kj}_{il}\ang a)\Omega_{il}
\end{equation}
in the algebra $\A\rti\A^\circ$. Comparing (\ref{comm1}) and (\ref{comm2}), we see that there is an $\A$-bimodule ${\rm map}~\tilde{\tau}:\tilde{\Gamma}_{v,\zeta}\rightarrow\A\rti\A^\circ$ such that $\tilde{\tau}(\theta_{kj})=\Omega_{kj}, k,j=1,{\dots},n$. Since $\tilde{\tau}(\theta)=C_{v,\zeta}+(Tr A)\varepsilon$ we derive 
$$
\tilde{\tau}(adb)=\tilde{\tau}(a(\theta b-b\theta))=a(C_{v,\zeta}b-b C_{v,\zeta}),a,b\in\A.
$$
Let $\tau$ be the restriction of $\tilde{\tau}$ to $\Gamma_{v,\zeta}$ and let $\rho$ be the algebra homomorphism of $\A$ into $L(\A)$ defined by 
$\rho(a)b=a\cdot b,b\in\A$. By the preceding we have proved the following

\mn
{\bf Proposition 4.} {\it The pair $(\rho,C_{v,\zeta})$ is an algebraic 
commutator representation of the bicovariant FODC $\Gamma_{v,\zeta}$ over $\A$ with $\B=L(\A)$.}

Let us specialize to the coordinate Hopf algebras $\A=\cO(G_q)$, where 
$G_q$ is one of the quantum groups $SL_q(n{+}1),O_q(n)$ or $Sp_q(2n)$; 
see [FRT] or [KS], Chapter 9. Since the Hopf algebra $\cO(G_q)$ is 
coquasitriangular ([KS], Theorem 10.9), the above considerations apply. 
In this case the algebraic commutator representation from Proposition 4 
is even faithful. The proof relies on the paper [HS] which was essentially based on results form [JL] and [J].

\mn
{\bf Theorem 5.} {\it Suppose that $q$ is transcendental. \\
(i) Let $G_q$ be  $SL_q(n+1),O_q(n)$ or $Sp_q(2n)$. Then the 
algebraic commutator representation $(\rho, C_{v,\zeta})$ of the FODC 
$\Gamma_{v,\zeta}$ is faithful.\\
(ii) Let $G_q=SL_q(n+1)$ or $G_q=Sp_q(2n)$. For any bicovariant 
finite dimensional FODC $\Gamma$ over $\cO(G_q)$ there exists 
a central element $C$ of the dual Hopf algebra $\cO(G_q)^\circ$ 
such that $(\rho,C)$ is a faithful algebraic commutator representation 
of $\Gamma$ with $\B=L(\A)$.}

\mn
{\bf Proof.} (i): Set $C:=C_{v,\zeta}$ and $\A:=\cO(G_q)$. Let 
$\{Y_1,{\dots},Y_m\}$ be a basis of the quantum tangent space 
$\tau_{v,\zeta}$ of the FODC $\Gamma_{v,\zeta}$. Then there are 
complex numbers $\alpha^s_{kj}$ such that 
$X_{kj}=\sum\alpha^s_{kj} Y_s$. Put 
$\Omega_s:=\sum_{k,j} \alpha^s_{kj}\Omega_{kj}$. We 
compute 
$$
\Delta(C)=\sum_{k,j} \zeta l^k_j\otimes \Omega_{kj}-(Tr A)\varepsilon\otimes\varepsilon
$$
and hence
\begin{equation}\label{ca}
C_{(1)}\langle C_{(2)},a\rangle-\langle C,a\rangle\varepsilon
=\sum_{k,j}\langle\Omega_{kj},a\rangle X_{kj}
=\sum_s\langle\Omega_s,a\rangle Y_s
\end{equation}
for $a\in\A$. In the proof of Theorem 4.1 in [HS], it was shown that the 
quantum tangent space $\tau_{v,\zeta}$ is the linear span of elements 
$C_{(1)}\langle C_{(2)}, a\rangle-\langle C,a\rangle\varepsilon$, 
where $a\in\A$. (The latter is equivalent to the assertion proved  in [HS] that the FODC $\Gamma_{v,\zeta}$ can be obtained by the method from [BM]. 
Note that in Theorem 4.1 in [HS] the quantum group $O_q(n)$ was excluded, but this part of the proof is valid for $O_q(n)$ as well.) Hence, by (\ref{ca}), the functionals $\langle\Omega_1,{\dots},\Omega_m\}$ of $\A^\circ$ are linearly independent. We have
$$
\tau(\sum_l a_ldb_l)=\sum_{l,k,j} a_l(X_{kj}\ang b_l)\Omega_{kj}=\sum_{l,s} a_l(Y_s\ang b_l)\Omega_s
$$
for $a_l,b_l\in\A$. Therefore, if $\tau(\sum_la_ldb_l)=0$, then it 
follows from the definition of the cross product algebra $\A\rti\A^\circ$ 
and the linear inpendence of $\{\Omega_1,{\dots},\Omega_m\}$ that 
$\sum_la_l(Y_s\ang b_l)=0$ for $s=1,{\dots},m$. This implies that 
$\sum_la_ldb_l=0$. Thus, $\tau$ is injective and $(\rho,C)$ is faithful.\\
(ii): By Theorem 4.1 in [HS], the quantum tangent space $\T_\Gamma$ of 
the FODC $\Gamma$ is the direct sum of quantum tangent spaces 
$\T_{v_j,\zeta_j},i=1,{\dots}, k$, of FODC $\Gamma_{v_j,\zeta_j}$.  Setting $C:=C_{v_1,\zeta_1}+{\dots}+C_{v_k,\zeta_k}$, the assertion follows by the repeating the reasoning of the proof of (i).
\hfill $\Box$

\mn
Next let us return to the case of general coquasitriangular Hopf 
algebras and consider Hilbert space commutator representations. 
For this we suppose that $\A$ is a coquasitriangular Hopf $\ast$-algebra 
and that the universal $r$-form $\br$ of $\A$ is real 
(that is, $\overline{\br(a\otimes b)}=\br(b^\ast\otimes a^\ast)$ for 
$a,b\in\A)$. Further, we assume that the matrix corepresentation 
$v=(v^j_k)$ is unitary (that is, $(v^j_k)^\ast=S(v^k_j)$ for 
$j,k=1,{\dots},n)$ and that the character $\zeta$ is hermitean (that is, $\zeta (a^\ast)=\overline{\zeta(a)}$ for $a\in\A$). Then the bicovariant FODC $\Gamma_{v,\zeta}$ is a $\ast$-calculus (see [KS], p. 520).

\mn
{\bf Lemma 6.} {\it Under the above assumptions, we have 
$(C_{v,\zeta})^\ast=C_{v,\zeta}$ in the $\ast$-algebra $\A^\circ$.}

\mn
{\bf Proof.} Using the well-known facts that $\br(S(a)\otimes S(b))= r(a\otimes b)$ and $S\circ \ast\circ S\circ\ast=\id$ we compute
\begin{align*}
\overline{A^j_k} &= \sum_s \overline{\br(S^2(v^j_s)\otimes v^s_k)}=\sum_s\br((v^s_k)^\ast\otimes (S^2(v^j_s))^\ast)\\
&=\sum_s \br(S(v^k_s)\otimes S^{-1}(S(v^j_s)^\ast))\\
&=\sum_s\br (S^2(v^k_s)\otimes S(v^j_s)^\ast)=\sum_s\br (S^2(v^k_s)\otimes v^s_j)=A^k_j.
\end{align*}
Since $(l^k_j)^\ast=l^j_k$ by formula (47) in [KS], p. 347, the assertion follows from (\ref{xkj}) and (\ref{cv}).\hfill $\Box$

We now proceed as in the preceding section. Let $\U_0$ be a 
$\ast$-subalgebra of $\A^\circ$ which contains all functionals 
$X_{kj},f^{kj}_{st},k,j,s,t=1,{\dots},n$, and let $\A\rti\U_0$
 denote the $\ast$-subalgebra of $\A\rti\A^\circ$ generated by $\A$ 
and $\U_0$. Retaining the above assumptions and notation, we have

\mn
{\bf Theorem 7.} {\it Suppose that $\pi_0$ is a $\ast$-representation of the $\ast$-algebra $\A\rti\U_0$. Let $\pi:=\pi_0\lceil\A$ and $F=\pi_0(C_{v,\zeta})$. Then the pair $(\pi,F)$ is a commutator representation of the bicovariant first order $\ast$-calculus $\Gamma_{v,\zeta}.$}

For $q\in\R$ the Hopf $\ast$-algebras $\cO(G_q)$ of the compact forms 
of quantum groups $G_q=SL_q(n+1),O_q(n), Sp_q(2n)$, are $CQG$ algebras 
and coquasitriangular with real universal $r$-form. Hence the preceding 
considerations apply to $\cO(G_q)$. Recall that the Haar state $h$ of 
the $CQG$ algebra $\cO(G_q)$ is faithful and that  $\pi_h$ denotes the 
Heisenberg representation of the cross product algebra 
$\cO(G_q)\rti \cO(G_q)^\circ$. Therefore, keeping the above assumption and 
combining Theorems 5 and 7 we obtain

\mn
{\bf Corollary 8.} {\it Suppose that $q\in \R$ is transcendental. \\
(i) Let $G_q$ be $SL_q(n{+}1),O_q(n)$ or $Sp_q(2n)$. Then the pair 
$(\pi_h,\pi_h(C_{v,\zeta}))$ is faithful commutator representation of the 
bicovariant $\ast$-FODC $\Gamma_{v,\zeta}$.\\
(ii) Let $G_q= SL_q(n{+}1)$ or $G_q=Sp_q(2n)$. For each finite 
dimensional bicovariant $\ast$-FODC $\Gamma$ there exists a central 
element $C\in\cO(G_q)^\circ$ such that $(\pi_h,\pi_h(C))$ is a 
faithful commutator representation of $\Gamma$.}

\bn
{\bf 4. Commutator Representations of FODC on Some Quantum Spaces}

There are a number of distinguished covariant $\ast$-calculi on 
quantum spaces which allow faithful commutator representations. In 
this final section we mention three such examples. 
We only state the corresponding formulas and omit the
straightforward verifications.

\sn
{\bf Example 1.} {\it Quantum disc, quantum complex plane}

Suppose that $\gamma \geq 0$ and $0<q<1$. Let $\X_{\gamma,q}$ denote 
the $\ast$-algebra with a single generator $z$ satisfying the relation 
$$
z^\ast z-q^2zz^\ast=\gamma(1-q^2).
$$
Note that $\X_{0,q}$ is the coordinate algebra of quantum complex 
plane, while $\X_{1,q}$ is the coordinate algebra of the quantum disc [KL]. 
On the $\ast$-algebra $\X_{\gamma,q}$ there is a distinguished 
$\ast$-calculus with bimodule structure given by 
$$
dz{\cdot} z=q^2 z dz, dz{\cdot} z^\ast=q^{-2}z^\ast dz, 
dz^\ast{\cdot} z=q^2 z dz^\ast, dz^\ast{\cdot} z^\ast=q^{-2}z^\ast dz^\ast.
$$
These simple relations have been found in [S2]. For $\gamma \neq 0$ the FODC $\Gamma$ has been extensively studied in [CHZ] and [SSV]. 

Define $C \in M_2(\X_{\gamma,q})$ and a homomorphism $\rho:\X_{\gamma,q}\to M_2(\X_{\gamma,q})$ by
\begin{align*}
C= (1-q^2)^{-1}\left(\begin{matrix} 0      &z\\ 
                                    z^\ast &0 
                     \end{matrix}\right),~ 
\rho(f) =\left(\begin{matrix} f&0\\ 
0&f\end{matrix}\right)
\end{align*}
for $f \in \X_{\gamma,q}$, we have
\begin{equation*}
[C,\rho(z)]=\left(\begin{matrix} 0&0\\ 1{-}z^\ast z                                   &0\end{matrix}\right),~ 
[C,\rho(z^\ast)]=\left(\begin{matrix} 0&1{-}z^\ast z\\ 
                                      0&0\end{matrix}\right).
\end{equation*}
Then the pair $(\rho,C)$ is a faithful algebraic commutator representation 
of the FODC with $\B=M_2(\X_{\gamma,q})$. 

Passing to a Hilbert space representation of 
the $\ast$-algebra $\X_{\gamma,q}$ we get a commutator representation 
of the FODC. Suppose now that $\gamma=1$. Then there is a well-known 
faithful $\ast$-representation $\pi$ of $\X_{\gamma,q}$ on a Hilbert space
with orthonormal basis $\{e_n,~n {\in} N\}$ such that
$$
\pi(z)e_n=\lambda_{n+1}e_{n+1},~\pi(z^\ast) e_n=\lambda_ne_{n-1},~
{\rm where}~ \lambda_n:=(1-q^{2n})^{1/2}.
$$
In particular, $\pi(1-z^\ast z)e_n=q^{2n+2}e_n$. This implies that 
for all $f,g \in \X_{\gamma,q}$ the operator $\pi(f)i[\pi(C),\pi(g)]$
is of trace class and so we have a 1-summable Fredholm module. 

Note that the calculus $\Gamma$ of $\X_{\gamma,q}$, $\gamma \neq 0$, 
is $\U_q(su(1,1))$-covariant. There are commutator representations which 
have the latter symmetry, but then the commutators $[F,\pi(g)]$ 
are not of trace class.

\bn
{\bf Example 2.} {\it Real quantum plane}

Suppose that $q$ is a complex number of modulus one. Let 
$\cO(\R^2_q)$ be the $\ast$-algebra with two hermitean 
generators $x{=}x^\ast, y{=}y^\ast$ and defining 
relation $xy{=}qyx$. We consider $\cO(\R^2_q)$ as $\ast$-subalgebra of 
the $\ast$-algebra 
${\hat \cO}(\R^2_q)$ with hermitean generators $x, y, y^{-1}$ and relations 
$xy=qyx$ and $yy^{-1}=y^{-1}y=1$.

There is a first order $\ast$-calculus $\Gamma$ of 
$\cO(\R^2_q)$ introduced in [PW], [WZ] and given by the relations
$$
xdx{=} q^{-2} dx{\cdot}x,~y dx{=}q^{-1} dx{\cdot}y+(q^{-2}{-}1)dy{\cdot}x,\\ 
y dx{=}q^{-1} dx{\cdot}y,~ y dy{=}q^{-2} dy{\cdot} y~.
$$
Let $\rho$ be as in Example 1 and define $C \in M_2({\hat \cO}(\R^2_q))$ by
\begin{align*}
&C=(q^2-1)^{-1}\left(\begin{matrix} q^2 x^2y^{-2} &0\\ 
0&y^{-2}\end{matrix}\right).
\end{align*}
Then we have
\begin{align*}
&[C,\rho(x)]=\left(\begin{matrix} q^2x^3y^{-2} &0\\ 
0&xy^{-2}\end{matrix}\right),
~[C,\rho(y)]=\left(\begin{matrix} x^2y^{-1} &0\\ 
0&0\end{matrix}\right).
\end{align*}
Then the pair $(\rho,C)$ is a faithful algebraic 
commutator representation of the FODC $\Gamma$ with 
$\B=M_2({\hat \cO}(\R^2_q))$.

Composed with a faithful $\ast$-representation $\pi$ of 
${\hat \cO}(\R^2_q)$, we obtain a faithful commutator representation of 
$\Gamma$. For instance, we may take the $\ast$-representation $\pi$ on 
the domain ${\cal D} :={\rm Lin}\{ e^{-\delta x^2{+} c x}; \delta \ge 0,
 c\in\C\}$ of the Hilbert space $L^2(\R)$ given by 
$x=e^{\alpha x}$ and $y=e^{\beta P}$, where $(Pf)(x)=-if^\prime(x)$, 
$q=e^{i\gamma}$, and $\alpha$ and $\beta$ are reals such that 
$\alpha \beta=\gamma$.

\bn
{\bf Example 3.} {\it Extended quantum plane} 

Suppose that $q\in\R$ and $q\neq 0$. Let ${\hat \cO} (\C^2_q)$ be the 
$\ast$-algebra with generators $x ,y$ and defining
relations
\begin{align*}
xy=qyx, xy^\ast=q^{-1}y^\ast x,
x^\ast x=xx^\ast, y^\ast y -y y^\ast=(q^{-2}-1)x^\ast x.
\end{align*}
The $\ast$-algebra ${\hat \cO} (\C^2_q)$ is 
the realification (see e.g. [KS], 10.2.7) of the coordinate algebra of the quantum plane $\C_q^2$ 
which is given by the relation $xy=qyx$. 
There is a first order $\ast$-calculus $\Gamma$ of 
$\hat{\cO} (\C^2_q)$ with bimodule structure described by the relations
\begin{align*}
dx{\cdot} x &= q^2 x{\cdot} dx,~ dx{\cdot} y = 
q y{\cdot} dx+(q^2-1)x{\cdot} dy,\\
dy{\cdot} y&=q^2 y{\cdot} dy,~ dy{\cdot} x=q x{\cdot} dy,\\
dx{\cdot} x^\ast&= q^{-2} x^\ast{\cdot} dx+(q^2-1)y^\ast{\cdot} dy,~
dx {\cdot}y^\ast=q^{-1} y^\ast dx,\\
dy{\cdot} x^\ast &= q^{-1} x^\ast{\cdot} dy,~ dy{\cdot} y^\ast=q^{-2} y^\ast dy.
\end{align*}
The first four relations are the same as in the preceding example. 
Note that the relations of this calculus are the ones described 
by Proposition 5 in [S3]. 

We first define a $\ast$-representation $\pi$ of 
$\hat{\cO} (\C^2_q)$. Let $\K$ be a fixed Hilbert space and let $\Hh=\oplus_{n=0}^\infty \Hh_n$, where
$\Hh_n:=\K$. For $\eta \in \K$ let $\eta_n$ denote the vector in $\Hh$ 
with $n$-th component $\eta$ and zero otherwise. 
Suppose that $N$ is a normal operator on the Hilbert space 
$\K$ with trivial kernel. Let $N=w|N|$ be the polar decomposition of $N$. 
Since ${\rm ker}N=\{0\}$ and $N$ is normal, $w$ is unitary. Put 
$\E= \cap_{k=0}^\infty \D(|N|^k)$. 
Then there is a $\ast$-representation $\pi$ of the $\ast$-algebra 
$\hat{\cO} (\R^2_q)$ on the domain 
$\D:={\rm Lin}\{\eta_n: \eta \in \E;n \geq 0\} $
given by
\begin{align*}
&\pi(x)\eta_n=q^{n+1}N\eta_n,
\pi(x^\ast)\eta_n=q^{n+1}N^\ast \eta_n,\\
&\pi(y)\eta_n=\lambda_n |N|\eta_{n-1},
\pi(y^\ast)\eta_n=\lambda_{n+1}|N|\eta_{n+1}.
\end{align*}
Suppose that $T$ and $S$ are operators on the dense domain $\E$ of the Hilbert space $\K$ 
such that $S$ is symmetric, $T(\E)\subseteq \E$, $T^\ast(\E)\subseteq \E$, $S(\E)\subseteq \E$ and satisfying the following conditions on the 
domain $\E$:
\begin{align*}
w^\ast Tw=qT,~ w^\ast Sw=q^2 S,~ T|N|=|N|T, S|N|=|N|S~.
\end{align*}
Let $F$ denote the symmetric operator on the domain $\D$ defined by
\begin{align*}
F\eta_n=\lambda_n T\eta_{n-1}+S\eta_n+\lambda_{n+1} T^\ast \eta_{n+1},~
\eta \in \K, n\in \N.
\end{align*}
Then the pair $(\pi,F)$ is a commutator representation of the 
$\ast$-calculus $\Gamma$.

\end{document}